\documentclass[12pt,reqno]{amsart}

\usepackage{amsmath,amsfonts,amssymb,amscd,amsthm,amsbsy,epsf}

 \newtheorem{thm}{Theorem}[section]

 \newtheorem{lem}[thm]{Lemma}
 
 \theoremstyle{definition}
 
 \theoremstyle{remark}
 
 \numberwithin{equation}{section}
 \newcommand{\C}{\mathbb{C}}

 \def\myskip{\noalign{\vskip6pt}}
\begin{document}

\title[Identities zeta function]
{Identities for the Riemann zeta function}

\author[Michael O. Rubinstein]
{Michael O.\ Rubinstein\\ \\
Pure Mathematics \\University of Waterloo\\200 University Ave W\\Waterloo, ON, N2L 3G1\\Canada}



%%% ----------------------------------------------------------------------

%\begin{abstract} We obtain several expansions for $\zeta(s)$ involving a
%sequence of polynomials in $s$, denoted in this paper by $\alpha_k(s)$. These
%polynomials can be regarded as a generalization of Stirling numbers of the
%first kind and our identities extend some series expansions for the zeta function that are
%known for integer values of $s$. The expansions also give a different approach
%to the analytic continuation of the Riemann zeta function.
%\end{abstract}

%%% ----------------------------------------------------------------------
\maketitle
%%% ----------------------------------------------------------------------

\section{Introduction}

In this paper, we obtain several expansions for $\zeta(s)$ involving a
sequence of polynomials in $s$, denoted by $\alpha_k(s)$. These
polynomials can be regarded as a generalization of Stirling numbers of the
first kind and our identities extend some series expansions for the zeta function that are
known for integer values of $s$. The expansions also give a different approach
to the analytic continuation of the Riemann zeta function.

The inspiration for our formulas comes from Kenter's short note in the May
1999 Monthly \cite{kenter} where he derives a formula for Euler's constant
$\gamma$ which can be regarded as the $s \to 0$ case of \eqref{eq:gamma2
intro}, after subtracting $1/s$ from both sides.

We start with Riemann's formula
\begin{equation}
    \label{eq:1}
    \begin{split}
    \Gamma (s) \zeta (s) & = \int_0^\infty \frac{x^{s-1}}{e^x-1}\,dx\ ,\qquad
    \Re s>1\ ,\\
    \myskip
    & = \int_0^\infty \frac{x^{s-1}e^{-x}}{1-e^{-x}}\,dx
    \end{split}
\end{equation}
and substitute $t=1-e^{-x}$ to get
\begin{equation}
    \label{eq:2}
    \int_0^1 (-\log (1-t))^{s-1} \frac{dt}t\ ,\qquad \Re s>1\ .
\end{equation}
Let
\begin{equation}
    g(t) = \frac{-\log (1-t)}t = \sum_0^\infty \frac{t^k}{k+1}\ ,\qquad |t|<1\ ,
\end{equation}
and consider the Taylor expansion
\begin{equation}
    \label{eq:3}
    g(t)^{s-1} = \sum_0^\infty \alpha_k (s) t^k\ ,\qquad |t|<1\ .
\end{equation}
%
%In this paper, we study the coefficients $\alpha_k(s)$ and
%apply the above series to the integral \eqref{eq:2}, and several variants, to obtain
%some identities for the Riemann zeta function.
%
We will derive the following recursion:
\begin{eqnarray}
    \label{eq:recur}
    \alpha_0(s)&=& 1 \ , \notag \\
    \alpha_1(s)&=& \frac{s-1}{2} \ , \notag \\
    \alpha_{k+1}(s)
    &=& \frac1{k(k+1)(k+2)} \sum_{j=1}^k
    \frac{\alpha_j(s)j(k+k^2+s(2k+2-j))}{(k-j+1)(k-j+2)} \ , \qquad k\geq 1\ , \notag\\
    &&
\end{eqnarray}
and prove the following theorem concerning $\alpha_k(s)$:
\begin{thm}

\label{thm:1}
For $k \geq 1$, $\alpha_k(s)/(s-1)$ is a polynomial in $s$ with
positive rational coefficients, and $\alpha_k(s)$ satisfies:
\begin{equation}
    |\alpha_k(s)| \leq c_s \frac{(1+\log(k+1))^{|s|+1}}{k+1}\ ,
\end{equation}
where
\begin{equation}
    c_s = \frac{|s-1|}{|s|+1} (|s|+2) 2^{|s|+1}\ .
\end{equation}

\end{thm}
We then have the following identities:
\begin{thm}
\label{thm:2}

\begin{equation}
    \label{eq:gamma2 intro}
    \Gamma (s) = \sum_0^\infty \frac{\alpha_k(s)}{s+k}\ .
\end{equation}
\begin{equation}
    \label{eq:main formula intro}
    \Gamma (s) \zeta (s) = \sum_0^\infty \frac{\alpha_k(s)}{s+k-1}\ .
\end{equation}
For positive integer $\lambda$:
\begin{eqnarray}
    \label{eq:3 intro}
    \Gamma(s) \zeta(s-\lambda)
    &=&
    \sum_{k=0}^\infty \alpha_k(s)
    \sum_{j=1}^\lambda
    (-1)^{\lambda+j} \frac{j! S(\lambda,j)}{s+k-j-1} \\
    &=&
    \sum_{k=0}^\infty \alpha_k(s) 
    \sum_{j=0}^{\lambda-1} E(\lambda,j) \frac{(\lambda-j-1)!}{(s+k-j-2)\ldots(s+k-\lambda-1)}\ , \notag \\
\end{eqnarray}
where $S$ and $E$ respectively denote the Stirling numbers of the second kind and the Eulerian numbers.
We also have:
\begin{equation}
    \label{eq:trigamma intro}
    \Gamma(s) \zeta(s+1)
    =
    \sum_{k=0}^\infty \alpha_k(s) \Psi_1(s+k)
    =
    \sum_{k=0}^\infty
    \sum_{n=0}^\infty
    \frac{\alpha_k(s)}{(s+k+n)^2}\ ,
\end{equation}
with $\Psi_1(s+k)$ the trigamma function.

\end{thm}

Equation \eqref{eq:main formula intro} is known in the case that $s$ is a positive integer
where the coefficients $\alpha_k(s)$ can be expressed in terms of the
Stirling numbers of the first kind. Jordan \cite[Sec. 68, (11), pg 194]{jordan} credits it to
Stirling. See also \cite{rs} and \cite{shen}. This is described in Section~\ref{sec:s1}.

\section{Recursions and bound for $\alpha_k(s)$}

To study the $\alpha_k(s)$'s in greater detail consider
\begin{equation}
    \frac{d}{dt} (g(t)^{s-1})
    = (s-1) g(t)^{s-2} \sum_0^\infty \frac{k+1}{k+2}\ t^k\ ,
\end{equation}
so 
\begin{equation}
    \label{eq:4}
    \frac{d}{dt} (g(t)^{s-1})g(t)
    = (s-1)g(t)^{s-1} \sum_0^\infty \frac{k+1}{k+2}\ t^k\ .
\end{equation}
On the other hand, differentiating \eqref{eq:3} term by term, the
left hand side above equals
\begin{equation}
    \label{eq:5} 
    \biggl( \sum_0^\infty (k+1)\alpha_{k+1}(s)t^k\biggr)
    \biggl( \sum_0^\infty \frac{t^k}{k+1}\biggr) \ .
\end{equation}
Equating coefficients of $t^k$ in \eqref{eq:4} and \eqref{eq:5}, we get
\begin{equation}
    \sum_{j=0}^k \frac{(j+1)\alpha_{j+1}(s)}{k-j+1}
    = (s-1) \sum_{j=0}^k \alpha_j (s)\frac{k-j+1}{k-j+2}\ .
\end{equation}
Using $\alpha_0(s)=1$, we write this as
\begin{equation}
    \label{eq:6} 
    \sum_{j=1}^{k+1} \frac{j-(s-1)(k-j+1)}{k-j+2} \alpha_j (s)
    =  (s-1) \frac{k+1}{k+2}\ ,\qquad k \geq 0.
\end{equation}
The first few $\alpha_j(s)$'s are listed below 
\begin{equation}
    \begin{split}
        &\alpha_0 (s) = 1\\
        &\alpha_1 (s) = \frac{s-1}2\\
        &\alpha_2 (s) = (s-1)\left( \frac18 s+\frac1{12}\right)\\
        &\alpha_3 (s) = (s-1)\left( \frac1{48}s^2 +\frac1{16}s +\frac1{24}\right)\\
        &\alpha_4 (s) = (s-1)\left(\frac1{384}s^3 +\frac7{384}s^2
        + \frac{23}{576}s + \frac{19}{720}\right)
    \end{split}
\end{equation}
We also observe, from \eqref{eq:6}, that we can write
\begin{equation}
    \label{eq:7}
    \alpha_{k+1}(s) = \frac{s-1}{k+2} -\frac1{k+1} \sum_{j=1}^k
    \frac{j-(s-1)(k-j+1)}{k-j+2} \alpha_j (s)\ .
\end{equation}
Furthermore, substituting $k-1$ for $k$ and moving the l.h.s. to the r.h.s.,
we have
\begin{equation}
    \label{eq:8}
    0 = \frac{s-1}{k+1} - \frac1k \sum_{j=1}^k \frac{j-(s-1)(k-j)}{k-j+1}
    \alpha_j (s)\ .
\end{equation}
Subtracting $(k+1)/(k+2)$ times \eqref{eq:8} from \eqref{eq:7}
we get
\begin{eqnarray}
    \alpha_{k+1}(s)
    &&= \frac1{k(k+1)(k+2)} \sum_{j=1}^k
    \frac{\alpha_j(s)j(k+k^2+s(2k+2-j))}{(k-j+1)(k-j+2)}\ , \qquad k\geq 1. \notag \\
    &&
\end{eqnarray}

This last form for the recursion governing the $\alpha_k$'s is useful in that
we easily obtain, inductively, the first part of Theorem~\ref{thm:1}, namely
that $\alpha_k(s)/(s-1)$, for $k \geq 1$, is a polynomial in $s$ with
positive rational coefficients.

Next, we determine a bound for $|\alpha_k(s)|$.
Let $M_s$ be the positive integer satisfying $|s|+1 \le M_s < |s|+2$.
From the positivity of the coefficients of $\alpha_k(s)/(s-1)$ we have
\begin{equation}
    \label{eq:inequality1}
    \left|\frac{\alpha_k(s)}{s-1}\right| < \frac{\alpha_k(M_s+1)}{M_s}\ ,
    \qquad k=1,2,\ldots.
\end{equation}
Here we are using the fact that $\alpha_k(s)/(s-1)$
is a polynomial in $s$ with positive coefficients and also
$M_s+1 >|s|$.
Rewriting \eqref{eq:inequality1}, and using $|s|+1 \le M_s$, we find
\begin{equation}
    \label{eq:9}
    |\alpha_k(s)| \le\frac{|s-1|}{M_s} \alpha_k (M_s +1) \leq
    \frac{|s-1|}{|s|+1} \alpha_k (M_s +1) \ .
\end{equation}
Now,

\begin{lem}
    \label{lem:2}
    $$
        |\alpha_k (M+1)| \le \frac{M2^{M-1}(1+\log (k+1))^{M-1}}{k+1}\ ,\qquad
        M= 1,2,3,\ldots
    $$
\end{lem}

\begin{proof}
$$
    \alpha_k (M+1)\text{ is the coefficient of $t^k$ in }
    g(t)^M = \left( 1+\frac{t}2 +\frac{t^2}3 +\cdots\right)^M \ .
$$
When $M=1$ the lemma is easily verified.
Now assume that the lemma has been verified for $M$ and consider the
$M+1$ case.
Writing $g(t)^{M+1} = g(t)^M g(t)$ we get
\begin{equation}
    \alpha_k (M+2) = \sum_{m=0}^k \frac1{m+1} \alpha_{k-m}(M+1)\ .
\end{equation}
By our inductive hypothesis,
\begin{align}
    |\alpha_k(M+2)| &\le \sum_{m=0}^k \frac1{m+1}
    \frac{M2^{M-1}(1+\log (k-m+1))^{M-1}}{k-m+1} \notag \\
    \myskip
    & = \frac{M2^{M-1}}{k+2} \sum_{m=0}^k (1+\log(k-m+1))^{M-1}
    \left(\frac1{m+1} + \frac1{k-m+1}\right) \notag \\
    \myskip
    & \le \frac{M}{k+2} 2^{M-1} (1+\log (k+1))^{M-1} \sum_{m=0}^k \frac2{m+1}\ .
\end{align}
But $2\sum_{m=0}^k 1/(m+1) \le2 (1+\log (k+1))$, and $M/(k+2) < (M+1)/(k+1)$, thus
the lemma is proven.
\end{proof}

Combining \eqref{eq:9} with Lemma~\ref{lem:2} and the assumption that $M_s < |s|+2$
yields the bound in Theorem~\ref{thm:1}.
\begin{equation}
    \label{eq:10}
    |\alpha_k(s)| \leq \frac{|s-1|}{|s|+1}\
    \frac{(|s|+2)2^{|s|+1} (1+\log(k+1))^{|s|+1}}{k+1}\ .
\end{equation}

We conclude, by using the above bound to show that the sum below converges,
that we may substitute \eqref{eq:3}
into \eqref{eq:2} and integrate the power series term by term to obtain
\begin{equation}
    \label{eq:main formula}
    \Gamma (s) \zeta (s) = \sum_0^\infty \frac{\alpha_k(s)}{s+k-1}\ .
\end{equation}
Away from its poles, the above is, in bounded sets, a uniformly convergent sum
of analytic functions. Thus, while we started in \eqref{eq:1} with $\Re s>1$,
formula~\eqref{eq:main formula} gives the meromorphic continuation of the left
hand side.

Similarly, applying the same change of variable to the integral
defining the Gamma function,
\begin{equation}
    \label{eq:gamma}
    \begin{split}
    & \int_0^\infty x^{s-1}e^{-x}\,dx \\
    &= \int_0^1 (-\log (1-t))^{s-1} dt\ ,\qquad \Re s>0\ ,
    \end{split}
\end{equation}
gives
\begin{equation}
    \label{eq:gamma2}
    \Gamma (s) = \sum_0^\infty \frac{\alpha_k(s)}{s+k}\ .
\end{equation}

\section{Connection to Stirling numbers of the first kind}
\label{sec:s1}

The coefficients $\alpha_k(s)$ defined by
\eqref{eq:3} are related to Stirling numbers of the first kind through the
series expansion, for integer $m \geq 0$,
\begin{equation}
    \label{eq:stirling}
    \log(1+t)^m = m! \sum_{n=m}^\infty \frac{s(n,m)}{n!}t^n.
\end{equation}

Thus, when $s$ is a positive integer,
\begin{equation}
    \label{eq:stirling2}
    \alpha_k(s) = (-1)^k \frac{(s-1)!}{(s+k-1)!} s(s+k-1,s-1)\ ,
\end{equation}
and equation \eqref{eq:main formula} becomes
\begin{equation}
    \label{eq:specialized}
     \zeta (s) = \sum_0^\infty (-1)^k \frac{s(s+k-1,s-1)}{(s+k-1)(s+k-1)!}\ .
\end{equation}
As mentioned in the introduction, the latter formula, valid for positive integer $s$,
is known
and \eqref{eq:main formula} may be regarded as a generalization of \eqref{eq:specialized}
to $s \in \C$.
The expression \eqref{eq:specialized} can be combined with identities that
relate the Stirling numbers to the harmonic numbers, and from this various
identities for zeta evaluated at positive integer values may be deduced. See
for example Section 4 of \cite{rs} or \cite{shen}.

\section{Connection to Stirling numbers of the second kind}

Other related expansions of $\zeta (s)$ are possible.
For example, instead of \eqref{eq:1}, consider
\begin{equation}
    \Gamma (s)\zeta (s-\lambda) =
    \int_0^\infty x^{s-1} 
    \sum_1^\infty n^\lambda e^{-nx}\,dx\ ,\qquad \Re s -\lambda >1\ .
\end{equation}
Substituting $t= 1-e^{-x}$, we can rewrite the above as
\begin{equation}
    \label{eq:variant}
    \int_0^1 t^{s-1} \left( \frac{-\log (1-t)}t\right)^{s-1}
    \sum_1^\infty n^\lambda (1-t)^{n-1}\,dt\ .
\end{equation}

When $\lambda$ is a positive integer, we can evaluate the sum over $n$
by repeatedly multiplying by $(1-t)$ and applying $-d/dt$
to the geometric series $\sum_1^\infty (1-t)^{n-1}=1/t$.
$$
    \vbox{\halign{\hfil$\displaystyle#$\hfil\qquad&\hfil$\displaystyle#$\hfil\cr
    \underline{\vphantom{\sum_1}\ \lambda\ }
    &\underline{\ \sum_1^\infty n^\lambda (1-t)^{n-1}\ }\cr
    \myskip
    1&1/t^2\cr
    \myskip
    2&2/t^3 -1/t^2\cr
    \myskip
    3&6/t^4 - 6/t^3 +1/t^2\cr
    \myskip
    4&24/t^5 - 36/t^4 +14/t^3 -1/t^2%\cr
    %\myskip
    %\vdots&\vdots\cr
    \cr}}
$$
Therefore,
\begin{align}
    \label{eq:a few variants}
    \Gamma (s)\zeta (s-1) & = \sum_0^\infty \frac{\alpha_k(s)}{s+k-2} \\
    \myskip
    \Gamma (s)\zeta (s-2) & = \sum_0^\infty \alpha_k(s)
    \left(\frac2{s+k-3} -\frac1{s+k-2}\right) \notag \\
    \myskip
    \Gamma (s)\zeta (s-3) & = \sum_0^\infty \alpha_k(s)
    \left(\frac6{s+k-4} -\frac6{s+k-3} +\frac1{s+k-2}\right)\notag\\
    \myskip
    \Gamma (s)\zeta (s-4) & =\sum_0^\infty \alpha_k(s) \left(\frac{24}{s+k-5}
    -\frac{36}{s+k-4} +\frac{14}{s+k-3} - \frac1{s+k-2}\right)\notag%\\
    %\myskip
    %\vdots\qquad & \quad \qquad \vdots \notag
\end{align}

Next, we derive the general form for the above expressions
and give a connection to Stirling numbers of the second kind:
For positive integer $\lambda$,
\begin{equation}
    \sum_1^\infty n^\lambda (1-t)^{n-1} =
    \sum_{j=1}^{\lambda}
    (-1)^{\lambda+j} j! S(\lambda,j)/t^{j+1}\ .
\end{equation}
We can
verify this formula inductively, multiplying by $(1-t)$, applying $-d/dt$, and using the
the recurrence relation for Stirling numbers of the second kind:
\begin{equation}
   \label{eq:S recursion}
   S(\lambda+1,j) = S(\lambda,j-1) + j S(\lambda,j)\ .
\end{equation}

Substituting into \eqref{eq:variant} and integrating gives the formula
\begin{equation}
    \label{eq:stirling 2}
    \Gamma(s) \zeta(s-\lambda)
    =
    \sum_{k=0}^\infty \alpha_k(s)
    \sum_{j=1}^\lambda
    (-1)^{\lambda+j} \frac{j! S(\lambda,j)}{s+k-j-1}
\end{equation}
of Theorem \ref{thm:2}.

We can also go in the opposite direction, taking, for example,
$\lambda=-1$ in \eqref{eq:variant}. Now,
\begin{equation}
    \sum_1^\infty n^{-1} (1-t)^{n-1} = \frac{\log{t}}{t-1}\ .
\end{equation}
Using
\begin{equation}
    \int_0^1 t^{s+k-1} \frac{\log{t}}{t-1} dt = \Psi_1(s+k)\ ,
\end{equation}
where $\Psi_1$ is the trigamma function, we get
\begin{equation}
    \label{eq:trigamma}
    \Gamma(s) \zeta(s+1)
    =
    \sum_{k=0}^\infty \alpha_k(s) \Psi_1(s+k)
    =
    \sum_{k=0}^\infty
    \sum_{n=0}^\infty
    \frac{\alpha_k(s)}{(s+k+n)^2}\ .
\end{equation}

\section{Connection to Eulerian numbers}

We can also expand the sum in \eqref{eq:variant}
in terms of Eulerian numbers. For positive integer $\lambda$,
\begin{equation}
    \label{eq:eulerian numbers}
    \sum_1^\infty n^\lambda (1-t)^{n-1} =
    t^{-\lambda-1} \sum_{j=0}^{\lambda-1}
     E(\lambda,j) (1-t)^{\lambda-j-1}\ ,
\end{equation}
where $E(\lambda,j)$ are the Eulerian numbers, satisfying the
recursion:
\begin{equation}
    E(\lambda+1,j) = (j+1) E(\lambda,j) + (\lambda+1-j) E(\lambda,j-1)\ .
\end{equation}
Formula \eqref{eq:eulerian numbers} can, again, be verified inductively,
multiplying by $(1-t)$, applying $-d/dt$, and manipulating slightly
% the manipulation is the following after multiplying and -d/dt, we
% get two sums from the product rule:
%     t^{-\lambda-2} \sum_{j=0}^{\lambda-1}
%     (\lambda+1) E(\lambda,j) (1-t)^{\lambda-j}\ ,
%     + t (\lambda-j) E(\lambda,j) (1-t)^{\lambda-j-1}\ ,
% we then subtract and add
% (\lambda-j) E(\lambda,j) (1-t)^{\lambda-j-1}\ ,
% and collect terms to get
%     t^{-\lambda-2} \sum_{j=0}^{\lambda-1}
%     (j+1) E(\lambda,j) (1-t)^{\lambda-j}\ ,
%     + (\lambda-j) E(\lambda,j) (1-t)^{\lambda-j-1}\ ,
% combine j term of the first sum and j-1 term of the second and apply the recursion.
before using the above recursion.

However,
\begin{equation}
    \label{eq:eulerian numbers b}
    \int_0^1 t^{s+k-\lambda-2}(1-t)^{\lambda-j-1} dt
    = \frac{\Gamma(s+k-\lambda-1)\Gamma(\lambda-j)}{\Gamma(s+k-j-1)}\ ,
\end{equation}
i.e. the Beta function, gives, after applying $\Gamma(z+1) = z \Gamma(z)$,
\begin{equation}
    \label{eq:eulerian2}
    \Gamma(s) \zeta(s-\lambda)
    =
    \sum_{k=0}^\infty \alpha_k(s) 
    \sum_{j=0}^{\lambda-1} E(\lambda,j) \frac{(\lambda-j-1)!}{(s+k-j-2)\ldots(s+k-\lambda-1)}\ .
\end{equation}

\section{Obtaining $\zeta(1-\lambda)$, $\lambda=1,2,3,\ldots$.}
\label{sec:special}

We can apply formula \eqref{eq:3 intro} to obtain expressions for $\zeta(1-\lambda)$,
with $\lambda$ a positive integer.
We first illustrate the technique for $\lambda=1$.

By Theorem~\ref{thm:1}, $\alpha_k(s)$, for $k \geq 1$,
is a polynomial in $s$ divisible by $s-1$. Thus,
$\alpha_k(1)=0$, $k \geq 1$. We also have $\alpha_0(1)=1$. Therefore, substituting
$s=1$ into the first equation in~\eqref{eq:a few variants}, all but the first {\it two} terms
vanish. The $k=0$ term gives $-1$, while the denominator of the $k=1$ term cancels the
$s-1$ factor of $\alpha_1(s)=(s-1)/2$, and the $k=1$ term equals $1/2$. We thus get
\begin{equation}
    \label{eq:zeta 0}
    \zeta(0) = -1 + 1/2 = -1/2\ .
\end{equation}

The general situation is handled via equation~\eqref{eq:3 intro}. As in the $\lambda=1$
case, we get two kinds of contributions- from the $k=0$ term, and from the terms
$k=j$, with $j=1,\ldots,\lambda$. The latter terms are the ones for which the denominator
$s+k-j-1$ cancels the factor $s-1$ of $\alpha_k(s)$.

The $k=0$ term produces, on simplifying,
\begin{equation}
    \label{eq:k=0 term}
    (-1)^\lambda
    \sum_{j=1}^\lambda
    (-1)^{j-1} (j-1)! S(\lambda,j)\ .
\end{equation}
But, the Stirling numbers of the second kind are defined by the expansion
\begin{equation}
    \label{eq:stirling2 defn}
    x^\lambda =
    \sum_{j=1}^\lambda
    S(\lambda,j) x (x-1) \ldots (x-j+1)\ , \qquad \lambda \geq 1.
\end{equation}
Thus, dividing by $x$, and setting $x=0$ shows that~\eqref{eq:k=0 term}
equals 0 if $\lambda>1$, and $-1$ if $\lambda=1$.
Therefore, the $k=0$ term only contributes when $\lambda=1$, and we denote its
contribution by $-\delta_\lambda$, equal to $0$ or $-1$ according to whether $\lambda>1$ or equals
$1$.

The terms $k=j$, with $j=1,2,\ldots,\lambda$, contribute
to~\eqref{eq:3 intro},
on canceling the zero and pole at $s=1$, the sum:
\begin{equation}
    \label{eq:other k}
    (-1)^\lambda \sum_{k=1}^\lambda \alpha_k(1)' (-1)^k k! S(\lambda,k)\ .
\end{equation}

Therefore , putting both contributions together:
\begin{equation}
    \label{eq:zeta 1-lambda}
    \zeta(1-\lambda) =
    -\delta_\lambda +
    (-1)^\lambda \sum_{k=1}^\lambda \alpha_k(1)' (-1)^k k! S(\lambda,k)\ .
\end{equation}

We thus need a formula for $\alpha_k(1)'$. We can
differentiate~\eqref{eq:7}, and use the fact from Theorem~\ref{thm:1} that $(s-1)\alpha_k(s)$
has a double order zero at $s=1$ when $k \geq 1$, to get
\begin{equation}
    \label{eq:alpha prime}
    \alpha_{k+1}(1)' = \frac{1}{k+2}
    -\frac1{k+1} \sum_{j=1}^k
        \frac{j}{k-j+2} \alpha_j(1)'\ , \qquad k \geq 0.
\end{equation}
This recursion allows us, in conjunction with~\eqref{eq:zeta 1-lambda}, to evaluate
$\zeta(1-\lambda)$ for any positive integer $\lambda$. 

Equation \eqref{eq:zeta 1-lambda}
can also be used
to derive Euler's famous formula involving the Bernoulli numbers.
To see the connection to Bernoulli numbers, we calculate the first few values 
of $\alpha_k(1)'$ from the above recursion, starting with $\alpha_1(1)'=1/2$, and list them
in Table~\ref{tab:alpha prime}.
\begin{table}[h!tb]
\centerline{\small
\begin{tabular}{|c|c|}
\hline
$k$ & $\alpha_k(1)'$ \cr \hline
0& 0 \cr
1& 1/2 \cr
2& 5/24 \cr
3& 1/8 \cr
4& 251/2880 \cr
5& 19/288 \cr
6& 19087/362880 \cr
7& 751/17280 \cr
8& 1070017/29030400 \cr
9& 2857/89600 \cr
10& 26842253/958003200 \cr
11& 434293/17418240 \cr
12& 703604254357/31384184832000 \cr
13& 8181904909/402361344000 \cr
14& 1166309819657/62768369664000 \cr
15& 5044289/295206912 \cr
16& 8092989203533249/512189896458240000 \cr
17& 5026792806787/342372925440000 \cr
18& 12600467236042756559/919636959090769920000 \cr
19& 69028763155644023/5377993912811520000 \cr
20& 8136836498467582599787/674400436666564608000000 \cr
\hline
\end{tabular}
}
\caption{Values of $\alpha_k(1)'$}
\label{tab:alpha prime}\end{table}

Googling some of the larger numerators in this table, for example
703604254357, immediately returns entries A002208 and A002657 from Sloane's Online
Encyclopedia of Integer Sequences \cite{Sl}. These entries deal with the numerators
of the Norlund numbers, and the following formula is stated:
`Numerator of integral of $x(x+1)...(x+n-1)$ from 0 to 1.'
Comparing a few of our numbers to the values of these integrals, we find that the
denominators in our case are off by a
factor of $k\ k!$ from the denominators of this formula.
We are thus led to surmise that:
\begin{lem}
    The following formula holds for $k>0$:
    \label{lem:formula alpha prime}
    \begin{equation}
        \alpha_k(1)' = \frac{1}{k\ k!}
        \int_{0}^1 (x)_k dx\ ,
    \end{equation}
    where
    \begin{equation}
        (x)_k = x(x+1)\ldots(x+k-1)\ .
    \end{equation}
\end{lem}
\begin{proof}
    We show that the r.h.s. in the lemma satisfies the same recursion as
$\alpha_k(1)'$. Replacing, in~\eqref{eq:alpha prime},
the $\alpha$'s by the r.h.s. of the lemma, rearranging and simplifying
slightly, we wish to show that
\begin{equation}
    \label{eq:goal lemma}
    \sum_{j=1}^{k+1} \frac{1}{k-j+2} \frac{1}{j!} \int_{0}^1 (x)_j dx\ =\
    \frac{k+1}{k+2}\ , \qquad k\geq 0.
\end{equation}
The l.h.s. above is the coefficient of $z^{k+2}$ in the product:
\begin{equation}
    \label{eq:product power series}
    \sum_{j=1}^\infty \frac{z^j}{j}\ \times \
    \sum_{j=1}^\infty \frac{z^j}{j!} \int_{0}^1 (x)_j dx\ .
\end{equation}
The first sum is the power series for $-\log(1-z)$, while the second sum equals the integral
from $x=0$ to 1 of
\begin{equation}
    \sum_{j=1}^\infty \frac{(x)_j}{j!} z^j\ ,
\end{equation}
which is the power series for $(1-z)^{-x}-1$. Integrating from 0 to 1,
we find that the product in~\eqref{eq:product power series} equals
\begin{equation}
    \label{eq:product2}
    - \log(1-z)
    \left(
        \frac{z}{(z-1)\log(1-z)} - 1
    \right)
    = \frac{z}{1-z} + \log(1-z),
\end{equation}
whose coefficient of $z^{k+2}$ is $1-1/(k+2) = (k+1)/(k+2)$.
\end{proof}

Substituting Lemma~\ref{lem:formula alpha prime} into \eqref{eq:zeta 1-lambda}
and rearranging summation and integration gives
\begin{equation}
    \label{eq:zeta 1-lambda 2}
    \zeta(1-\lambda) =
    -\delta_\lambda +
    (-1)^\lambda \int_0^1 \sum_{k=1}^\lambda (-1)^k \frac{S(\lambda,k)}{k} (x)_k dx\ .
\end{equation}
Moving the $(-1)^k$ into the $(x)_k$ and changing variables, $u=-x$, yields
\begin{equation}
    \label{eq:zeta 1-lambda 3}
    -\delta_\lambda +
    (-1)^\lambda \int_{-1}^0 \sum_{k=1}^\lambda \frac{S(\lambda,k)}{k} u(u-1)\ldots(u-k+1) du\ .
\end{equation}
Next, apply the recursion \eqref{eq:S recursion}
to split the above sum over $k$ into two sums. The second of these sums equals
\begin{equation}
    \sum_{k=1}^{\lambda-1} S(\lambda-1,k) u(u-1)\ldots(u-k+1)\ ,
\end{equation}
which, by \eqref{eq:stirling2 defn} is $u^{\lambda-1}$, if $\lambda>1$.
If $\lambda=1$ it equals 0. Integrating from $-1$ to $0$, it contributes
to~\eqref{eq:zeta 1-lambda 3}:
\begin{equation}
    \begin{cases}
        -1/\lambda\, \qquad \text{if $\lambda >1$,} \\
        0, \qquad \text{if $\lambda=1$.} \notag
    \end{cases}
\end{equation}
which we can write as
\begin{equation}
        -1/\lambda +\delta_\lambda\ .
\end{equation}
Therefore,~\eqref{eq:zeta 1-lambda 3} has been simplified to
\begin{equation}
    \label{eq:zeta 1-lambda 4}
    \zeta(1-\lambda) = -1/\lambda +
    (-1)^\lambda \int_{-1}^0 \sum_{k=1}^\lambda \frac{S(\lambda-1,k-1)}{k} u(u-1)\ldots(u-k+1) du\ .
\end{equation}
However, the sum above
can be related to sums of powers, and hence to Bernoulli polynomials.
Assume for now that $u$ is a positive integer.
First, we use the operator $\Delta f(m) = f(m+1)-f(m)$, applied to
$m(m-1)\ldots(m-k+1)$ to get rid of the numerator above:
\begin{equation}
    (m+1)\ldots(m-k+2) - m(m-1)\ldots(m-k+1) = k m (m-1) \ldots (m-k+2).
\end{equation}
Dividing by $k$, summing over $m=0$ to $u-1$, and telescoping yields
the well known identity
\begin{equation}
    \frac{u(u-1)\ldots(u-k+1)}{k} = \sum_{m=0}^{u-1} m (m-1) \ldots (m-k+2), \qquad k \geq 1.
\end{equation}
(The $m=0$ term is needed if $k=1$.)
Therefore, the sum in \eqref{eq:zeta 1-lambda 4} equals
\begin{equation}
    \sum_{k=1}^\lambda S(\lambda-1,k-1) \sum_{m=0}^{u-1} m (m-1) \ldots (m-k+2).
\end{equation}
Rearranging the two sums, and using
\begin{equation}
    \sum_{k=1}^\lambda S(\lambda-1,k-1) m (m-1) \ldots (m-k+2) = m^{\lambda-1}\ ,
\end{equation}
the integrand in~\eqref{eq:zeta 1-lambda 4} becomes, for positive integer $u$,
\begin{equation}
    \sum_{m=0}^{u-1} m^{\lambda-1}.
\end{equation}
But, sums of powers can be expressed in terms of the Bernoulli polynomials and numbers:
\begin{equation}
    \label{eq:bernoulli identity}
    \sum_{m=0}^{\lambda-1} m^{\lambda-1} = \frac{B_\lambda(u) - B_\lambda }{\lambda}\ , \qquad \lambda \geq 1,
\end{equation}
giving a formula for the integrand which is valid for all $u$, and not just positive integer $u$,
since both the integrand in~\eqref{eq:zeta 1-lambda 4}
and the above are polynomials in $u$ agreeing on infinitely many values.
We thus have
\begin{equation}
    \label{eq:zeta 1-lambda 5}
    \zeta(1-\lambda) = -1/\lambda +
    (-1)^\lambda \int_{-1}^0
    \frac{B_\lambda(u) - B_\lambda }{\lambda} du\ , \qquad \lambda \geq 1.
\end{equation}
Furthermore,
\begin{equation}
    \label{eq:bernoulli identity 2}
    \int_{-1}^0 B_\lambda(u) du = (-1)^\lambda
\end{equation}
which can be obtained by substituting $t=u+1$, using the difference equation for the Bernoulli polynomials,
\begin{equation}
    \label{eq:bernoulli diff}
    \frac{B_\lambda(u+1)-B_\lambda(u)}{\lambda} = u^\lambda, \quad \lambda \geq 1,
\end{equation}
and applying the third of the defining properties of the Bernoulli polynomials:
\begin{eqnarray}
    \notag
    B_0(t) &=& 1 \\
    \notag
    B_k'(t) &=& k B_{k-1}(t), \quad k \geq 1 \\
    \int_0^1 B_k(t) dt &=& 0 , \quad k \geq 1.
\end{eqnarray}
Plugging \eqref{eq:bernoulli identity} \eqref{eq:bernoulli identity 2} into
~\eqref{eq:zeta 1-lambda 5} therefore gives Euler's formula:
\begin{equation}
    \label{eq:zeta 1-lambda Euler}
    \zeta(1-\lambda) = (-1)^{\lambda+1} \frac{B_\lambda}{\lambda}\ , \qquad \lambda \geq 1.
\end{equation}

\section{Further properties of $\alpha_k$}

Because $B_\lambda=0$ for odd $\lambda>1$,
Euler's formula gives, as is well known, $\zeta(-m)=0$ for positive even integer $m$.
Hence, the left hand side of \eqref{eq:main formula} has poles
at $s=1,0,-1,-3,-5,-7,\ldots$.
Therefore, to cancel the poles of the right hand side at $s=-2,-4,-6,\ldots$,
$\alpha_k(s)$ must be divisible by $s+k-1$ when $k=3,5,7,\ldots$.

Similarly, from the first equation in \eqref{eq:a few variants},
$\alpha_k(s)$ is divisible by $s+k-2$ when
$k=1,3,5,7,\ldots$, and from the second formula, $2\alpha_{k+1}(s) -\alpha_k(s)$ is
divisible by $s+k-2$ when $k=0,2,4,6,\ldots$.

By comparing the residue of~\eqref{eq:gamma2} at $s=-k$, $k \geq 0$, we get
\begin{equation}
    \label{eq:alpha minus k}
    \alpha_k(-k) = (-1)^k/k!.
\end{equation}
Next, by considering the residue of \eqref{eq:main formula} at the poles $s=-1,-3,-5,\ldots$,
we have
\begin{equation}
    \alpha_{2m+2} (-2m-1) = \frac{B_{2m+2}}{(2m+2)!} \ ,\qquad m \geq 0,
\end{equation}
and from first formula in \eqref{eq:a few variants},
\begin{equation}
    \alpha_{2m+2}(-2m) = - (2m+1) \frac{B_{2m+2}}{(2m+2)!}\ ,\qquad
    m \geq 0.
\end{equation}

% ------------------------------------------------------------------------
\bibliographystyle{amsplain}

% ------------------------------------------------------------------------

\end{document}